\documentclass[11pt]{article}
\newcommand{\dated}{\mbox{} \hfill {\small [{\tt \today}]}}
\usepackage{graphicx,float}
\restylefloat{figure}
\newtheorem{theorem}{Theorem}
\title{Why prove things?}
\author{\textit{Volker Runde}}
\date{}
\begin{document}
\maketitle
Once upon a time, there was a prince who was educated by private tutors. One day, the math tutor set out to explain the Pythagorean theorem to his royal student. The prince wouldn't believe it. So, the teacher proved the theorem, but the prince was not convinced. The teacher presented another proof of the theorem, and then yet another, but the prince would still shake his head in disbelief. Desperate, the teacher exclaimed: ``Your royal highness, I give you my word of honor that this theorem is true!'' The prince's face lit up: ``Why didn't you say so right away?!''
\par
Wouldn't that be wonderful? A simple word of honor from the teacher, and the student accepts the theorem as true\ldots
\par 
Of course, it would be awful. Who makes sure that the teacher can be trusted? Where did he get his knowledge from? Did he rely on another person's word of honor? Was the person from whom the teacher learned the theorem trustworthy? Where did that person get his/her knowledge from? Did that person, too, trust someone elses's word of honor? The longer the chain of words of honor gets, the shakier the theorem starts to look. It can't go on indefinitely: someone must have established the truth of the theorem some other way. My guess is: that someone \emph{proved} it.
\par 
Why are mathematicians so obsessed with proofs? The simple answer is: because they are obsessed with the truth. A proof is a procedure which, by applying certain rules, establishes an assertion as true. Proofs do not only occur in mathematics. In a criminal trial, for instance, the prosecution tries to \emph{prove} that the defendant is guilty. Of course, the rules according to which a proof is carried out depend very much on the context: it is one thing to prove in court that Joe Smith stole his neighbor's hubcaps and another one to give a proof that there are infinitely many prime numbers. But in the end all proofs serve one purpose: to get to the truth.
\par 
Here is a joke: \textit{A mathematician and a physicist are asked to check whether all odd numbers greater than one are prime. The mathematician says: ``Three's a prime, five's a prime, seven's a prime, but nine isn't. Therefore it's false.'' The physicist says: ``Three's a prime, five's a prime, seven's a prime, nine isn't---but eleven and thirteen are prime again. So, five out of six experiments support the hypothesis, and it's true!''} We laugh at the physicist. How can he simply dismiss a counterexample? It's not as silly as it seems. Experimental data rarely fits theoretical predictions perfectly, and scientists are used to a certain amount of data that is somewhat out of line. What makes the physicist in the joke look foolish is that he treats a mathematical problem like an experimental one: he applies rules of proof that are valid in one area to another area where they don't work.
\par 
Instead of musing further on the nature of mathematical proof, let's try and do one.
\par 
Consider a chessboard consisting of 64 squares:
\begin{figure}[H]
\begin{center}
\includegraphics[width=0.33\linewidth,clip]{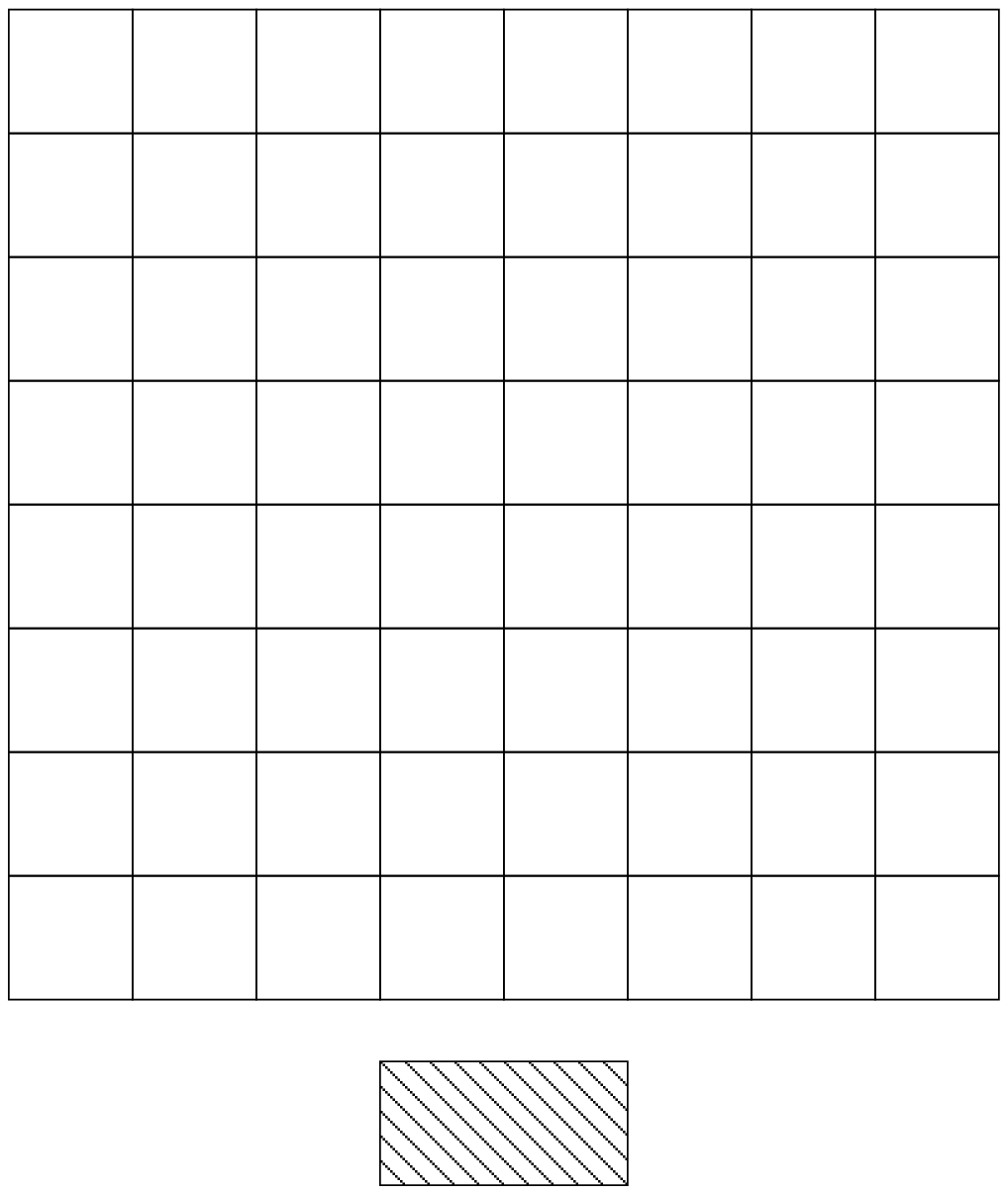} 
\end{center}
\end{figure}
Then take rectangular tiles as shown below the board: each of them covers precisely two adjacent squares on the chessboard. It's obvious that you can cover the entire board with such tiles without any two of them overlapping. That's straightforward, so why do we need proof here? Not yet\ldots
\par
To make things slightly more complicated, take a pair of scissors to the chessboard and cut away the squares in the upper left and in the lower right corner. This is how it will look like:
\begin{figure}[H]
\begin{center}
\includegraphics[width=0.33\linewidth,clip]{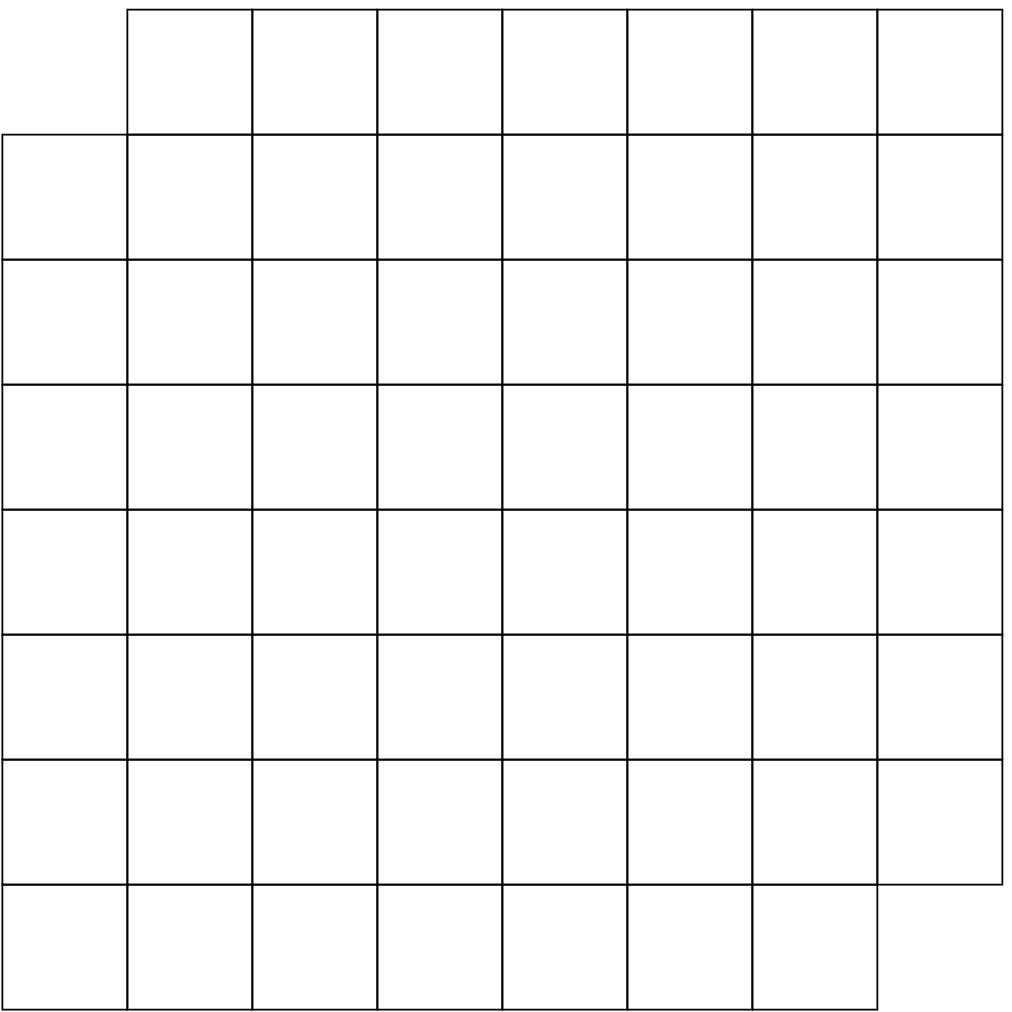} 
\end{center}
\end{figure}
Now, try to cover this altered chessboard with the tiles without any two of them overlapping\ldots
\par 
If you really try this (preferably with a chessboard drawn on a piece of paper\ldots), you'll soon find out that---to say the least---it's not easy, and maybe the nagging suspicion will set in that it's not even possible---but why? 
\par 
There is, of course, the method of brute force to find out. There are only finitely many ways to place the tiles on the chessboard, and if we try them all and see that in no case the area covered by them is precisely the altered chessboard, then we are done. There are two problems with this approach: firstly, we need to determine \emph{every} possible way to arrange the tiles on the chessboard, and secondly, even if we do, the number of possible tile arrangements may be far too large for us to check them all. So, goodbye to brute force\ldots
\par 
So, if brute force fails us, what can we do? Remember, we are dealing with a chessboard, and a chessboard not only consists of 64 squares in an eight by eight pattern---the squares alter in color; 32 are white, and 32 are black:
\begin{figure}[H]
\begin{center}
\includegraphics[width=0.33\linewidth,clip]{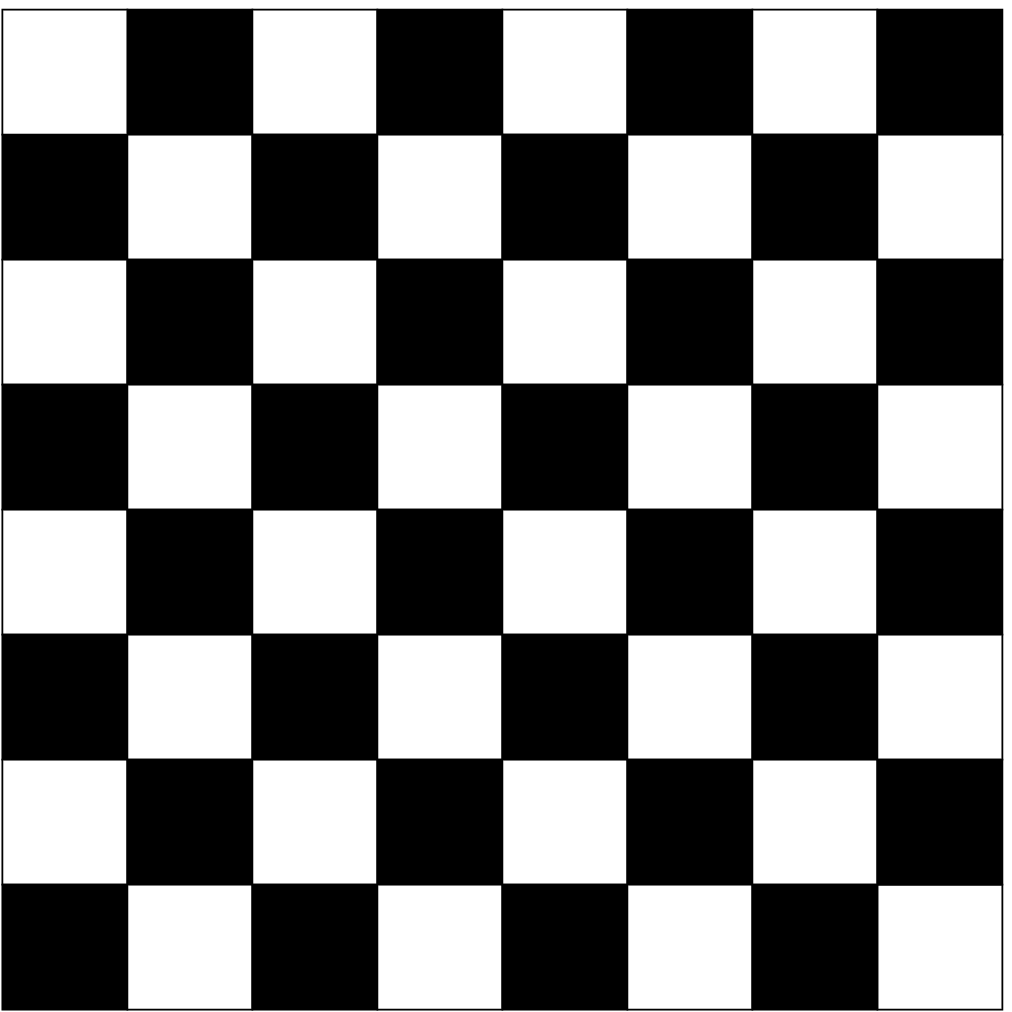} 
\end{center}
\end{figure}
Each tile covers precisely two adjacent squares on the board, and two adjacent squares on a chessboard are always different in color; so each tile covers one white square and one black square. Consequently, any arrangement of tiles on the chessboard must cover the same number of white and of black squares. But now, check the altered chessboard: 
\begin{figure}[H]
\begin{center}
\includegraphics[width=0.33\linewidth,clip]{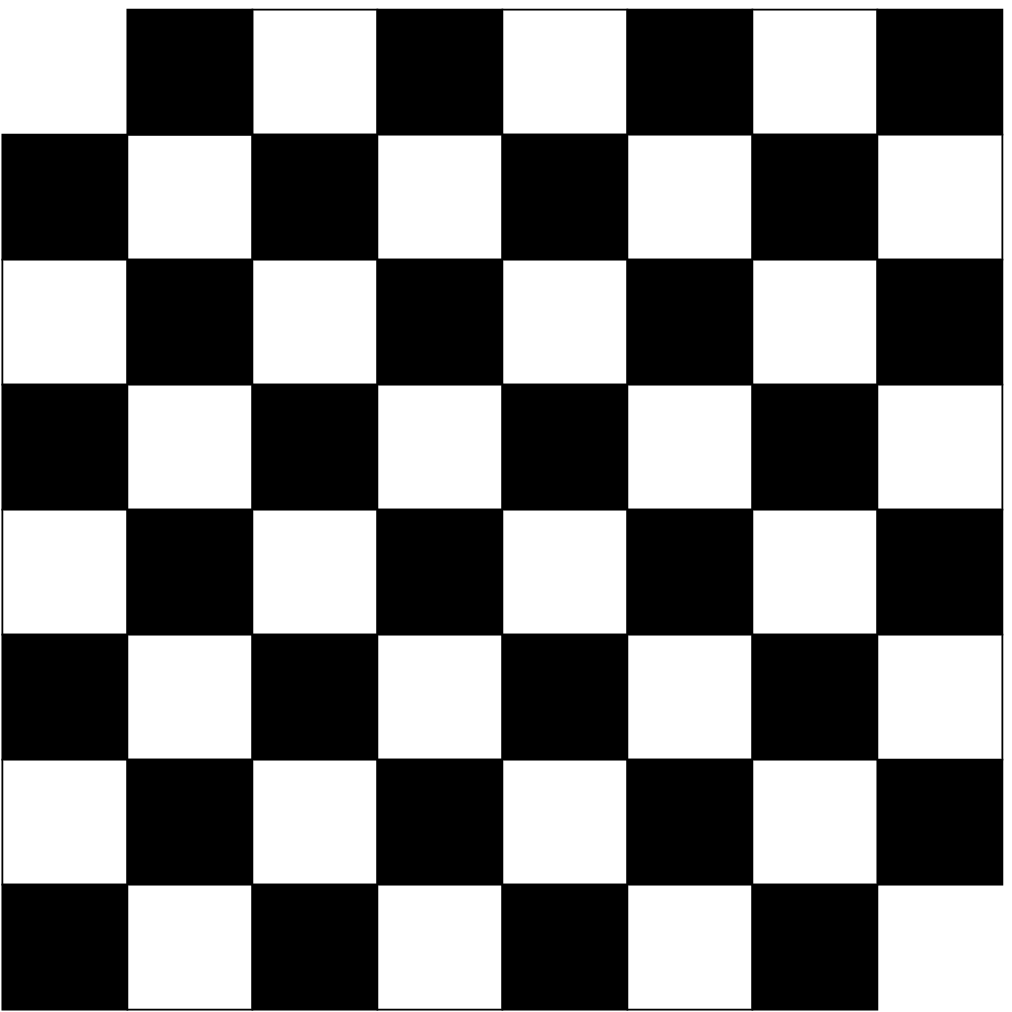} 
\end{center}
\end{figure}
We removed two white squares, so the altered board has 30 white squares, but 32 black ones. Therefore, it is impossible to cover it with tiles---we \emph{proved} it.
\par 
Does this smell a bit like black magic? Maybe, at the bottom of your heart, you prefer the brute force approach: it's the harder one, but it's still doable, and maybe you just don't want to believe that the tiling problem is unsolvable unless you've tried \emph{every} possibility.
\par 
There are situations, however, where a brute force approach to truth is not only inconvenient, but impossible. Have a look at the following mathematical theorem:
\begin{theorem} \label{thm1}
Every integer greater than one is a product of prime numbers.
\end{theorem}
\par 
Is it true? And if so, how do we \emph{prove} it?
\par 
Let's start with checking a few numbers: $2$ is prime (and thus a product of prime numbers), so is $3$, $4=2 \cdot 2$, $5$ is prime again, $6= 2 \cdot 3$, $7$ is prime, $8 = 2 \cdot 2 \cdot 2$, $9 = 3 \cdot 3$, and $10 = 2 \cdot 5$. So, the theorem is true for all integers greater than $2$ and less than or equal to $10$. That's comforting to know, but what about integers greater than $10$? Well, $11$ is prime, $12 = 2 \cdot 2 \cdot 3$, $13$ is prime, $14 = 2 \cdot 7$, $15= 3 \cdot 5$, $16 = 2 \cdot 2 \cdot 2 \cdot 2$, \ldots I stop here because it's useless to continue like this. There are infinitely many positive integers, and no matter how many of them we can write as a product of prime numbers, there will always remain infinitely many left for which we haven't shown it yet. Is $10^{10^{10^{10}}}+1$ a product of prime numbers? That number is awfully large. Even with the help of powerful computers, it might literally take an eternity to find the prime numbers whose product it is (if they exist\ldots). And if we have shown that the theorem holds true for every integer up to $10^{10^{10^{10}}}+1$, we still don't know about $10^{10^{10^{10^{10}}}}+1$.
\par 
Brute force leads nowhere here. Checking the theorem for certain examples might give you a feeling for it---but it doesn't help to establish its truth for \emph{all} integers greater than one.
\par 
Is the theorem possibly wrong? What would that mean? If not every integer greater than one is a product of prime numbers, than there must be at least one integer $a_0$  which is \emph{not} a product of prime numbers. Maybe, there is another integer $a_1$ with $1 < a_1 < a_0$ which is also not a product of prime numbers; if so replace $a_0$ by $a_1$. If there is an integer $a_2$ with $1 < a_2 < a_1$ which is not a product of prime numbers, replace $a_1$ by $a_2$. And so on\ldots There are only finitely many numbers between $2$ and $a_0$, and so, after a finite number of steps, we hit rock bottom and wind up with an integer $a > 1$ with the following properties: (a) $a$ is not a product of prime numbers, and (b) it is the smallest integer with that property, i.e., every integer greater than one and less than $a$ is a product of prime numbers.
\par
Let's think about this (hypothetical) number $a$. It exists if the theorem is false. What can we say about it? It can't be prime because then it would be a product (with just one factor) of prime numbers. So, $a$ isn't prime, i.e., $a=bc$ with neither $b$ nor $c$ being $a$ or $1$. This, in turn, means that $1 < b,c < a$. By property (b) of $a$, the numbers $b$ and $c$ are thus products of prime numbers, i.e., there are prime numbers $p_1, \ldots, p_n, q_1, \ldots, q_m$ such that $b = p_1 \cdots p_n$ and $c = q_1 \cdots q_m$. But then
\[
  a = bc = p_1 \cdots p_nq_1 \cdots q_m
\]
holds, and $a$ is product of prime numbers, which contradicts (a).
\par 
We assumed that the theorem was wrong, and---based on that assumption---obtained an integer $a$ that is not a product of prime numbers only to see later that this was not possible. The only way out of this dilemma is that our assumption was wrong: the theorem is true! (And we now know that $10^{10^{10^{10}}}+1$ is a product of prime numbers without having to find them\ldots)
\par 
The strategy we used to prove Theorem \ref{thm1} is called \emph{indirect proof}. We can't show something directly, so we assume it's wrong and (hopefully) arrive at a contradiction.
\par 
Let's try another (indirect) proof:
\begin{theorem} \label{thm2}
There are infinitely many prime numbers.
\end{theorem}
\par 
Is this believable? There is no easy formula to calculate the $n$th prime number, and after putting down the first few prime numbers, it gets harder and harder to come up with the next prime. So, is the theorem wrong and do we simple run out of prime numbers after a while? 
\par 
Assume this is so: there are only finitely many prime numbers, say $p_1, \ldots, p_n$. Set $a := p_1 \cdots p_n + 1$. By Theorem \ref{thm1}, $a$ is a product of prime numbers. In particular, there are a prime number $q$ and a non-negative integer $b$ with $a = qb$. Since $p_1, \ldots, p_n$ are all the prime numbers there are, $q$ must be one of them. Let $c$ be the product of all those $p_j$ that aren't $q$, so that $a = q c+1$. We then obtain
\[
  0 = a - a = qc+1 - qb = q(c-b)+1,
\]
and thus $q(c-b) = -1$. This, however, is impossible because $c-b$ is a non-zero integer and $q \geq 2$.
\par 
We have thus again reached a contradiction, and Theorem \ref{thm2} is proven.
\par
The proof of Theorem \ref{thm2} isn't as straightforward as the one of Theorem \ref{thm1}. Why did we define $a$ the way we did? The answer is simply that it works this way, and it's been working for over two thousand years: Theorem \ref{thm2} was first stated (and proven) in Euclid's \textit{Elements}, which appeared around 300 B.C. As the American mathematician Saunders Mac Lane once said: ``Mathematics rests on proof---and proof is eternal.''
\\[\bigskipamount]\dated
\vfill
\begin{tabbing}
\textit{Address}: \= Department of Mathematical and Statistical Sciences \\
\> University of Alberta \\
\> Edmonton, Alberta \\
\> Canada T6G 2G1 \\[\medskipamount]
\textit{E-mail}: \> \texttt{vrunde@ualberta.ca} \\[\medskipamount]
\textit{URL}: \> \texttt{http://www.math.ualberta.ca/$^\sim$runde/}   
\end{tabbing} 
\end{document}